\documentclass{article}

\newcommand{\R}{I\!\! R}

\begin{document}

Tsemo Aristide

Centre de Recherche Mathematiques

Universite de Montreal

Case Postale 6128, Succursale Centre Ville

Montreal Quebec, H3C 3J7

\bigskip
\bigskip

\centerline{\bf Connective structures for principal gerbes.}

\bigskip
\bigskip

 \centerline{\bf Introduction.}

\medskip

Let $N$ be a manifold, $H$ a Lie group and $P$ a $H$-principal
bundle defined over $N$, a $P$-gerbe defined over $N$, is a gerbed
defined over $N$ bounded by the sheaf of automorphisms of $P$. In
this paper we define the fundamental notions of the differential
geometry of $P$-gerbes, that is the notions of connective
structure curving holonomy and characteristic classes. The
interest of such definitions is to provide a geometric action in
string theory. In classical physics, the variational functional of
the evolution of a particle in a phase space is a function of the
holonomy of a gauge connection. In string theory, the action is
given by the holonomy of a Deligne connective structure, when the
gauge group is the circle. It is natural to provide the definition
of non abelian holonomy in order to describe the action when the
gauge group is non commutative.

\bigskip

{\bf Acknowledgements.}

The authors want to thank Pierre Deligne for helpful corrections,
and Dusa McDuff and Johan Dupont for helpful discussions.

\medskip

{\bf  Notations.}

\medskip

Let $U_{i_1},...,U_{i_p}$ be open subsets  of a manifolds $N$, and
$C$ a presheaf defined on $N$. We will denote by $U_{i_1..i_p}$
the intersection of $U_{i_1}$,...,$U_{i_p}$. If $e_{i_1}$ is an
object of $C(U_{i_1})$, ${e_{i_1}}^{i_2...i_p}$ will be the
restriction of $e_{i_1}$ to $U_{i_1...i_p}$. For a map
$h:e\rightarrow e'$ between two objects of $C(U_{i_1..i_p})$, we
denote by $h^{i_{p+1..i_n}}$ the restriction of $h$ to a morphism
between $e^{i_{p+1}...i_n}\rightarrow {e'}^{i_{p+1}...i_n}$.

\medskip

 {\bf Definition.}

Let $N$   be a  manifold, $H$ a Lie group and $P\rightarrow N$ a
$H$-principal bundle defined on $N$. A {\bf $P$-gerbe} is a gerbe
bounded by the sheaf of automorphisms of $P$. More precisely it is
defined as follows:

To each open subset $U$ of $N$ we associate a category $C_P(U)$.
The group of automorphisms of an object of $C_P(U)$ is the group
of automorphisms of the restriction of $P$ to $U$. We suppose that
the following conditions are satisfied:

Gluing conditions for objects.

Let $U$ be an open subset of $N$,  $(U_i)_{i\in I}$ an open cover
of $U$ and $e_i$ an object of $C_P(U_i)$. Suppose given an arrow
$u_{ij}:e^i_j\rightarrow e^j_i$ between the respective
restrictions of $e_j$ and $e_i$ to $U_i\cap U_j$ such that
${u_{i_1i_2}}^{i_3}{u_{i_2i_3}}^{i_1}={u_{i_1i_3}}^{i_2}$. Then
there exists an object $e_U$ of $C_P(U)$ whose restriction to
$U_i$ is $e_i$.

\medskip

Gluing conditions for arrows.

Let $e$ and $e'$ be two objects of $C_P(U)$. The correspondence
defined on the category of open subsets of $U$ by $V\rightarrow
Hom(e_{\mid V},e'_{\mid V})$ is a sheaf of sets, where $e_{\mid
V}$ and $e'_{\mid V}$ are the respective restrictions of $e$ and
$e'$ to $V$.

\medskip

We suppose that there exists an open cover of $N$ $(U_i)_{i\in I}$
such that the category $C_P(U_i)$ is not empty, and objects of
$C_P(U_i)$ are isomorphic.

\medskip

An example of $P$-gerbe is defined as follows: Let $G$ be a Lie
group and $H$ a closed normal subgroup of $G$. The quotient $G/H$
is a Lie group. We suppose that  the projection $G\rightarrow G/H$
has local sections. Consider a $G/H$-bundle $p_{G/H}:P\rightarrow
N$, defined on the manifold $N$. That is a locally trivial bundle
whose transition functions is defined by the trivialization
$(U_i,u_{ij})$, $u_{ij}:U_i\cap U_j\rightarrow G$ defined by the
coordinate changes:

$$
U_i\cap U_j\times G/H\longrightarrow U_i\cap U_j\times G/H
$$

$$
(x,y)\longrightarrow (x,yu_{ij}(x))
$$

The functions $u_{ij}$ verify the following property:
${u_{i_3i_1}}^{i_3}(x){u_{i_1i_2}}^{i_3}(x){u_{i_2i_3}}^{i_1}(x)$
are in the center of  $H$ to insure the $G/H$-bundle to be
well-defined  as the $H$-bundle $p_H$  whose transition functions
are defined by:

$$
U_i\cap U_j\times H\longrightarrow U_i\cap U_j\times H
$$

$$
(x,y)\longrightarrow (x,{u_{ij}}^{-1}(x)y{u_{ij}}(x))
$$

Let ${\cal H}$ be the Lie algebra of $H$ and $Ad$ the adjoint
representation. We can define the locally trivial ${\cal
H}$-bundle $p_{\cal H}$ over $N$ whose transition functions are
defined by:

$$
U_i\cap U_j\times {\cal H}\longrightarrow U_i\cap U_j\times {\cal
H}
$$

$$
(x,y)\longrightarrow (x,Ad({u_{ij}}^{-1})(x))(y))
$$

\medskip

{\bf Proposition.}

{\it  Let  $U$ be an open subset of $N$, we denote by  $C_H(U)$
the category of $G$-principal bundles whose quotient by $H$ is the
restriction of $p_{G/H}$ to $U$. A morphism between a pair of
objects $e$ and $e'$ of $C_H(U)$ is a morphism of $G$-bundles
which cover the identity of the restriction of $p_{G/H}$ to $U$.
The correspondence  defined on the category of open subsets of $N$
by $U\rightarrow C_H(U)$ is a gerbe bounded by the sheaf of
automorphisms of $p_H$.}

\medskip

{\bf Proof.}

Gluing property for arrows:

Let $U$ be an open subset of $N$, and $(U_i)_{i\in I}$ an open
cover of $U$. Consider an object $e_i$ of $C_H(U_i)$, and a map
$u_{ij}:e^i_j\rightarrow e^j_i$ such that
${u_{i_1i_2}}^{i_3}{u_{i_2i_3}}^{i_1}={u_{i_1i_3}}^{i_2}$ The
definition of bundle implies the existence of a  $G$-bundle $e$
whose restriction to $U_i$ is $e_i$. Since the quotient of $e_i$
by $H$ is the restriction of $p_{G/H}$ to $U_i$, we deduce that
the quotient of $e$ by $H$ is the restriction of $p_{G/H}$ to $U$.

Gluing condition of arrows.

Let $e$ and $e'$ be a pair of objects of $C_H(U)$, the
correspondence defined on the category of open subsets of $U$ by
$V\rightarrow Hom(e_{\mid V},e'_{\mid V})$ is a sheaf of sets,
since it is the sheaf of morphisms between two bundles. The
bundles $e_{\mid V}$ and $e'_{\mid V}$ are the respective
restrictions of $e$ and $e'$ to $V$.

\medskip

Consider a trivialization $(U_i,u_{ij})$ of $p_{G/H}$. The bundle
$U_i\times G$ is an element of $C_H(U_i)$, thus $C_H(U_i)$ is not
empty, and for each object $e$ and $e'$ of $C_H(U)$, the
restrictions of $e$ and $e'$ to $U_i\cap U$ are isomorphic to the
trivial bundle $U_i\cap U\times G$ by an isomorphism whose
projection to $U_i\cap U\times G/H$ is the identity.

Let $e$, be an object of $C_H(e)$, and $f$ an automorphism of $e$,
The restriction $f_i$ of $f$ to the restriction of $e$ to $U_i\cap
U$ is an automorphism of the trivial bundle $U_i\cap U\times G$
which projects to the identity on $U_i\cap U\times G/H$. We deduce
that $f_i$ is defined by a map $f':U_i\cap U_j\rightarrow H$. On
$U_i\cap U_j\cap U$, we have $f_j={u_{ij}}^{-1}f_i{u_{ij}}$. This
implies that $f$ is a section of $p_{H}$ $\bullet$

\bigskip

{\bf The classifying cocycle of a principal gerbe.}

\medskip

Let $(U_i)_{i\in I}$ be an open cover of $N$ such that the
category $C_P(U_i)$ is not empty and the objects of $C_P(U_i)$ are
isomorphic. Consider for each $i$, an object $e_i$ of $C_P(U_i)$,
and arrow $u_{ij}:e^i_j\rightarrow e^j_i$. We can defined the
automorphism of $e^{i_1i_2}_{i_3}$:
$c_{i_1i_2i_3}={u_{i_3i_2}}^{i_1}{u_{i_2i_1}}^{i_3}{u_{i_1i_3}}^{i_2}$.

\medskip

{\bf Proposition.}

{\it The family of maps $c_{i_1i_2i_3}$ is a non commutative Cech
$2$-cocycle.}

\medskip

{\bf Proof.}

Let
${c'_{i_1i_2i_3}}^{i_4}={u_{i_4i_3}}^{i_1i_2}{c_{i_1i_2i_3}}^{i_4}{u_{i_3i_4}}^{i_1i_2}$
On
$U_{i_1i_2i_3i_4}$, we have

$$
c'_{i_1i_2i_3}c_{i_1i_3i_4}=c_{i_2i_3i_4}c_{i_1i_2i_4}
$$

$\bullet$

{\bf Connective structure on $H$-gerbes.}

\medskip

{\bf Definition.}

  Consider a gerbe $C_P$ defined on a manifold $N$ whose band is $L$, the sheaf of
  automorphisms
  of the principal bundle $P\rightarrow N$.
  {\bf A connective structure} on $C_P$, is a correspondence which
associates to each object $e_U$ of $C_P(U)$ an affine space
$Co(e_U)$, called the torsor of connections, which is a subset of
the set of ${p_{{\cal H}}}_{\mid U}$-valued $1$-forms defined on
$U$, where ${p_{{\cal H}}}_{\mid U}$ is the restriction of
$p_{{\cal H}}$ to $U$. The following properties are supposed to be
satisfied by this assignment:

(i)- The correspondence $e_U\rightarrow Co(e_U)$ is functorial
with respect to restrictions  to smaller subsets.

(ii)- For every isomorphism $h:e_U\rightarrow e'_U$ between
objects of $C_P(U)$, there exists an isomorphism of torsors
$h^*:Co(e_U)\rightarrow Co(e'_U)$ compatible with the composition
of morphisms of $C_P(U)$, and the restrictions to smaller subsets.

 (iii)- For each morphism $g$ of the object
$e_U$ of $C_P(U)$, and $\nabla_{e_U}$ a connection of $Co(e_U)$,

$$
g^*\nabla_{e_U}=Ad(g^{-1})(\nabla_{e_U})+g^{-1}dg
$$

For each open subset $U$, we define $C'_P(U)$ to be the category
whose objects are pair of objects $(e_U,\nabla_{e_U})$, where
$\nabla_{e_U}$ is an element of $Co(e_U)$. A morphism
$f:(e_U,\nabla_{e_U})\rightarrow (e'_U,\nabla_{e'_U})$ is
$\nabla_{e'_U}-u^*(\nabla_{e_U})$ where $u:e_U\rightarrow e'_U$ is
a morphism of $C(U)$. We suppose the correspondence $U\rightarrow
C'_P(U)$ to be a gerbe $\bullet$

\bigskip

{\bf The classifying cocycle of a connective structure.}

\medskip

 Let $(U_i, h_{ij})_{i\in I}$ a trivialization of
$p_{H}$ such that $C_P(U_i)$ is not empty, and $e_i$ an object of
$C_P(U_i)$ and $u_{ij}$ a morphism between $e^i_j$ and $e^j_i$.
Consider an element $\alpha_i$ of $Co(e_i)$. We define
$c_{i_1i_2i_3}$ to be
${u_{i_3i_1}}^{i_2}{u_{i_1i_2}}^{i_3}{u_{i_2i_3}}^{i_1}$. On
$U_i\cap U_j$, we can define the ${\cal H}$-valued form
$\alpha_{ij}=\alpha^j_i-{u_{ij}}^*(\alpha^i_j)$.

The Cech boundary of the $p_{{\cal H}}$ $1$-cocycle $\alpha_{ij}$
is:

$$
{u_{i_1i_2}}^*(\alpha_{i_2i_3})-\alpha_{i_1i_3}+\alpha_{i_1i_2}
$$

$$
{u_{i_1i_3}}^*(\alpha_{i_3}-{{c_{i_1i_2i_3}}}^*(\alpha_{i_3}))
$$

$$
=Ad({h_{i_1i_3}}^{-1})(\alpha_{i_3}-Ad({c_{i_1i_2i_3}}^{-1})(\alpha_{i_3})+{c_{i_1i_2i_3}}^{-1}d({c_{i_1i_2i_3}}))
$$

We have used the fact that on the trivialization $U_i\cap
U_j\times {\cal H}$, let $\alpha$ and $\alpha'$ be two elements of
$Co(e_i)$, ${u_{ij}}^*(\alpha-\alpha')$ is transformed in
$Ad({h_{ij}}^{-1})(\alpha-\alpha')$ by the transition functions of
$p_{{\cal H}}$, since $\alpha-\alpha'$ is an element of the vector
space of the affine space $Co(e_i)$.

\medskip

\medskip

{\bf Example.}

\medskip

Consider a normal subgroup $H$ of a Lie group $G$, and a
$G/H$-bundle $p_{G/H}$ over the manifold $N$. We have defined a
gerbe $C_H$ at page 2. We define the connective structure $Co$ on
$C_H$ as follows: for each open subset $U$ of $N$, and an object
$e_U$ of $C_H(U)$, $Co(e_U)$ is the set of $1$-forms
$\theta':U\rightarrow C({\cal H})$ where $C({\cal H})$ is the
center of ${\cal H}$. This definition is natural since if the
center of $H$ is trivial, then the gerbe $C_H$ is trivial. The
characteristic classes defined below are also trivial.

\medskip

{\bf Definition.}

{\bf A curving} of a connective structure $Co$ is a correspondence

$$
D(e_U,):Co(e_U)\rightarrow D(e_U)
$$

 where $D(e_U)$ is an affine space
whose underlying vector space is a set of $p_{\cal H}$ valued
$2$-forms which satisfies the following property:

(i)- For each morphism $h:e'_U\rightarrow e_U$,
$(D(e_U,\nabla))=D(e'_U,h^*\nabla)$.

(ii)- If $\alpha$ is a ${p_{{\cal H}}}_{\mid U}$  $1$-form on $U$
such that $\nabla+\alpha$ is an element of $Co(e_U)$, then

$$
D(e_U,\nabla+\alpha)=D(e_U,\nabla)+d\alpha+\alpha\wedge\alpha
$$

The assignment $e_U\rightarrow D(e_U,\nabla)$ is compatible with
the restrictions to smaller subsets.

\medskip

 {\bf
Characteristic classes.}

\medskip

{\bf Definition.}

A polynomial function of degree $l$ $F:{\cal H}^l\rightarrow {\R}$
is
 said to be invariant if and only if for every $h\in { H}$
$F(Ad(h))=F$.

Let $C_P$ be a gerbe bounded by $P$ endowed with the connective
structure $C'_P$, consider $[\Omega]$ the cohomology class of the
classifying cocycle $\Omega$ of $C'_P$ identified with a DeRham
$3$-form using the Cech-DeRham isomorphism. For every invariant
polynomial $P$ of degree $l$ we can define the $3l$-form
$P(\Omega)$ by $P\circ \wedge^lL$. The cohomology classes of the
forms $P(\Omega)$ are the characteristic classes of the curving.

\bigskip

{\bf Holonomy of non abelian gerbes.}

\medskip

Let $C_P$ be a principal gerbe defined over a manifold $N$ endowed
with a connective structure $Co$ and a curving $Cur$. Let
$l:N_2\rightarrow N$ be a differentiable map whose domain is the
compact surface $N_2$. We can pull-back the gerbe $C_P$, $Co$ and
$Cur$ to $N_2$ by using $l$. Let $(U_i)_{i\in I}$ be an open
covering of $N_2$ such that $(l^*C_P)(U_i)$ is not empty. Consider
an object $e_i$ of $(l^*(C_P)(U_i))$, $\nabla_i$ an element of
$l^*(Co)(U_i)$, and $L_i$ the curving of $\nabla_i$. Since $N_2$
is a surface, $L_i$ is exact. We can set $L_i=d(L'_i)$. On
$U_{i_1i_2}$ we have
$L_{i_2}-L_{i_1}=d(\nabla_{i_1}-u_{i_1i_2}^*\nabla_{i_2})+(\nabla_{i_1}-u_{i_1i_2}^*\nabla_{i_2})
\wedge(\nabla_{i_1}-u_{i_1i_2}^*\nabla_{i_2})$, the form
$(\nabla_{i_1}-u_{i_1i_2}^*\nabla_{i_2})\wedge
(\nabla_{i_1}-u_{i_1i_2}^*\nabla_{i_2})=d(L'_{i_1i_2})$ This
implies that
$L'_{i_2}-L'_{i_1}=\nabla_{i_1}-u_{i_1i_2}^*\nabla_{i_2}+L'_{i_1i_2}+
d(L"_{i_1i_2})$. This implies that the Cech boundary
$\delta(h_{i_1i_2})$ of
$h_{i_1i_2}=\nabla_{i_1}-u_{i_1i_2}^*\nabla_{i_2}+L'_{i_1i_2}$ is
a $2$-chain of closed forms. We set
$\delta(h_{i_1i_2})=d(C_{i_1i_2i_3})$ The chain
$C_{i_1i_2i_3}+\delta(L"_{i_ii_2})$ is a $2$-chain of constant
${\cal H}$-functions. Since $N_2$ is a surface, we can find a
cover such that this chain is a cocycle. It suffices to find an
open cover $(U_i)_{i\in I}$ such that  $U_{i_1i_2i_3i_4}$ is empty
and $U_i$ is a $1$-Eilenberg-Mclane space. Thus using the
Cech-DeRham isomorphism, we  identify this chain to a $l^*(p_{\cal
H})$ $2$-form $H$. We define

$$
Hol(N_2,C_P,Co)=exp(\int_{N_2}H)
$$

\medskip

\centerline{\bf Reference.}

\medskip

J.L Brylinski, Loops spaces, Characteristic Classes and Geometric
Quantization, Progr. Math. 107, Birkhauser, 1993.

\end{document}